\documentclass[11pt,a4paper]{amsart}
\usepackage{amsmath,amssymb, amsbsy}
\usepackage{color,psfrag}
\usepackage[dvips]{graphicx}
\usepackage{color}
\usepackage{mathtools}
\usepackage{hyperref}
\theoremstyle{plain}

\newtheorem{theorem}{Theorem}[section]

\newtheorem{corollary}{Corollary}[section]
\newtheorem{remark}{Remark}[section]

\newcommand{\hn}{\mathbb{H}^{N}}

\newcommand{\dv}{\: {\rm d}v_{\hn}}

\newcommand{\dvg}{\: {\rm d}v_{g}}
\newcommand{\dr}{\: {\rm d}r}

\newcommand{\dx}{\: {\rm d}x}

\numberwithin{equation}{section} \allowdisplaybreaks

\usepackage[text={6in,8.6in},centering]{geometry}
\begin{document}
	
	\title[Improved Poincar\'e-Hardy and CKN Inequalities]{Improved Poincar\'e-Hardy inequalities  \\on certain subspaces of the Sobolev space}

	\author[Debdip Ganguly]{Debdip Ganguly}
	\address{\hbox{\parbox{5.7in}{\medskip\noindent{Department of Mathematics,\\
					Indian Institute of Technology Delhi,\\
					IIT Campus, Hauz Khas, Delhi,\\
					New Delhi 110016, India. \\[3pt]
					\em{E-mail address: }{\tt 
						debdipmath@gmail.com}}}}}
					
	\author[Prasun Roychowdhury]{Prasun Roychowdhury}
	\address{\hbox{\parbox{5.7in}{\medskip\noindent{Department of Mathematics,\\
					Indian Institute of Science Education and Research,\\
					Dr. Homi Bhabha Road, Pashan,\\
					Pune 411008, India. \\[3pt]
					\em{E-mail address: }{\tt prasunroychowdhury1994@gmail.com}}}}}
	
	\date{\today}

	\maketitle
	
	\begin{abstract}
		We prove an improved version of Poincar\'e-Hardy inequality in suitable subspaces of the Sobolev space on the hyperbolic space via Bessel pairs.  As a consequence, we obtain a new Hardy type inequality with an improved constant (than the usual Hardy constant). Furthermore, we derive a new kind of improved Caffarelli-Kohn-Nirenberg inequality on the hyperbolic space.
	\end{abstract}
	

\section{Introduction}
This paper aims to study an improved version of the family of Poincar\'e inequalities in a suitable subspace of the Sobolev spaces on the hyperbolic space $\hn$ and by an improvement, we mean adding a Hardy type remainder term for the Poincar\'e inequality.   
Cartan-Hadamard manifolds are known to admit the following Hardy inequality (See \cite{Carron}) :
\begin{align}\label{hardy-manifold}
	\int_M|\nabla_g u|^2\dvg\geq\frac{(N-2)^2}{4}\int_M\frac{u^2}{\varrho(x,o)^2}\dvg \text{ for all }u\in C_c^\infty(M),
\end{align}
where $\varrho$ denotes the geodesic distance and $\nabla_g,\dvg$ represents the Riemannian gradient and volume element in $(M,g).$ Subsequently several authors studied and improved Hardy-type inequalities on manifolds,  (we refer few of them, e.g., \cite{Dambrosio, pinch, Kombe1, YSK}). Moreover, it is well known that the constant $\frac{(N-2)^2}{4}$ in \eqref{hardy-manifold} is optimal and never achieved in the corresponding  Sobolev space.

\medspace

It is worth mentioning that Cartan-Hadamard manifolds whose sectional curvatures are bounded above by a strictly negative constant admit Poincar\'e inequality. Let  $\hn$ be the real hyperbolic space, a simplest example of a Cartan-Hadamard manifold having a constant 
negative sectional curvature. Let $\lambda_{1}(\hn)$ with $N\geq 2$ denotes the bottom of the $L^2-$spectrum of the Laplace-Beltrami operator $-\Delta_{\hn}$ on the hyperbolic space, we write as 
\begin{equation*}\label{poincare}
	\lambda_{1}(\hn)= \inf_{u \in C_{c}^{\infty}(\hn)  \setminus \{ 0 \}} \dfrac{\int_{\hn} |\nabla_{\hn} u|^2 \  \emph{\rm d}v_{\hn}}{\int_{\hn} u^2 \  \emph{\rm d}v_{\hn}}=\left(\frac{N-1}{2} \right)^2 \,.
\end{equation*}
Moreover, the infimum is never achieved and this suggests an improvement might possible.  In recent years, there have been many efforts to improve the Poincar\'e inequality, namely by adding the remainder term of Hardy type and its generalization to Cartan-Hadamard manifolds. Poincar\'e-Hardy inequality has been studied thoroughly by many authors, we refer to (\cite{AK, EDG, BGGP, AKR, VHN}) and the references quoted therein.

\medspace

Recently in \cite{FLLM}, the author studied the family of Poincar\'e-Hardy inequalities on Cartan-Hadamard manifolds using Bessel pairs and their results provided an improvement of those already known in the literature. Their main novelty lies in proving identities by substituting suitable Bessel pairs in the non-flat setting. In the case of the Euclidean space, Frank-Seiringer in \cite{Frank} exploited the nonlinear and non-local version of the {\it ground state representation \rm} to derive  the following identity : 
\begin{align*}
	\int_{\mathbb{R}^N}V(x)|\nabla u|^2 \dx = 	\int_{\mathbb{R}^N} W(x) |u|^2 \dx + \int_{\mathbb{R}^N} V(x) |f|^2 \bigg|\nabla \bigg(\frac{u}{f}\bigg)\bigg|^2 \dx \quad \forall \ u \in C_c^{\infty}(\mathbb{R}^N),
\end{align*}
where $f$ is a positive weak solution of the Laplace equation
\begin{align*}
	-\text{div}(V(x)\nabla f(x))=W(x)f(x).
\end{align*}

\medspace

In pursuit of the above study, Ghoussoub-Moradifam in \cite{GM} established Hardy inequality with radial weights by introducing the notion of the Bessel pair, 
namely $C^1$ functions $(V,W)$ on $(0,R)$ is said to be a Bessel pair if the ordinary differential equation
\begin{align*}
	(Vy^\prime)^\prime+Wy \, = \, 0,
\end{align*}
admits a positive solution on $(0,R)$. If $(r^{N-1}V,r^{N-1}W)$ is a Bessel pair subject to the constraints that $\int_{0}^{R} \frac{1}{r^{N-1} V(r)} \, {\rm d}r = \infty$ and $\int_{0}^{R} r^{N-1} V(r) \, {\rm d}r < \infty$ where $r=|x|$ and $0<R\leq \infty$ is the radius of the ball $B$, then on ball $B$ in $\mathbb{R}^N$ the following inequality holds true for some positive constant $C >0:$
\begin{equation*}\label{1eq}
	\int_{B} V(x) |\nabla u|^2 \, {\rm d}x \geq C \int_{B} W(x) \, |u|^2 \, {\rm d}x \quad \forall \ u \in C_c^{\infty}(B).
\end{equation*}


\medskip 

In context to Hardy inequalities, with particular choices of $(V, W)$ the results in \cite{GM} simplified and improved several known results concerning Hardy inequalities 
and their $L^2$- improvements in the literature. Recently exploiting the notion of Bessel pairs the authors in \cite{LamTP, LLZ} established improved Hardy identities and inequalities in different spaces like on upper-half paces and Homogeneous 
groups. Motivated from the results of \cite{GM}, as mentioned before,  an attempt has been made in \cite{FLLM} to obtain improved Hardy inequalities on Cartan-Hadamard manifolds $(M,g)$ using  Bessel pairs. It was shown that for $u\in C_0^\infty(B_R(o)\setminus \varrho^{-1}\{0\})$ there holds
\begin{align*}
	\int_{B_R(o)}V(\varrho(x))|\nabla_g u|^2{\rm d}v_g&\geq\int_{B_R(o)}W(\varrho(x))|u|^2 {\rm d}v_g\\&+\int_{B_R(o)}V(\varrho(x))\varphi^2(\varrho(x))\bigg|\nabla_g\bigg(\frac{u}{\varphi(\varrho(x))}\bigg)\bigg|_g^2{\rm d}v_g,
\end{align*}
where $\varrho(x)= \mbox{dist}(x,o)$ with a fixed pole $o\in M$. Also note that in the above inequality,  $V$ and $W$  are positive $C^1$-functions on $(0,R)$ such that $(r^{N-1}V,r^{N-1}W)$ forms a Bessel pair with non increasing positive solution $\varphi$. These results also hold for the operator if we replace it with its radial part (refer \cite[Theorem 2.1]{FLLM}). Moreover, authors in  \cite{FLLM} established several interesting Hardy identities and inequalities which in particular generalise many well-known Hardy inequalities 
on Cartan-Hadamard manifolds, for example ({\cite{AK, EG, BGGP, EGR, Dambrosio, Kombe1, AKR, VHN, prc-20}}).

\medspace

Drawing primary motivation from \cite{FLLM}, recently a variant of  Hardy-Poincar\'e identities have been established in \cite[Theorem~6.1]{EGR-21}, which reads as : let $N \geq 2$, then for all $0 \leq \lambda \leq \lambda_{1}(\hn)$  and for all $u \in C_c^{\infty}(\hn \setminus \{ x_0\})$ there holds 
	\begin{align}\label{improved-poinc-lambda}
		\int_{\hn} |\nabla_{\hn}u|^2 \, {\rm d}v_{\hn} &= \lambda \int_{\hn} u^2 \, {\rm d}v_{\hn} + h^2_{N}(\lambda) \int_{\hn} \frac{u^2}{r^2} \, {\rm d}v_{\hn}  \\
		& + \left[\frac{(N-2)^2}{4} - h^2_{N}(\lambda)  \right] \int_{\hn} \frac{u^2}{\sinh^2 r} \, {\rm d}v_{\hn} \notag\\&+ \gamma_N(\lambda) h_N(\lambda) \int_{\hn} \frac{r \coth r - 1}{r^2} \,u^2 \, {\rm d}v_{\hn}   \notag \\
		&+ \int_{\hn}(\Psi_{\lambda}(r))^2\bigg|\nabla_{\hn}\bigg(\frac{u}{\Psi_{\lambda}(r)}\bigg)\bigg|^2\dv,\notag
	\end{align}
	where $\gamma_{N}(\lambda):=\sqrt{(N-1)^2-4\lambda}$, $h_N(\lambda):=\frac{\gamma_{N}(\lambda)+1}{2},\Psi_{\lambda}(r) := r^{-\frac{N-2}{2}} \left(\frac{\sinh r}{r} \right)^{-\frac{N-1+\gamma_N(\lambda)}{2}}$ and $r:={\rm d}(x, x_0)$ is the geodesic distance from a fixed pole $x_0 \in \hn$. It is worth mentioning that \eqref{improved-poinc-lambda}, first appeared as an inequality in \cite{BGGP}.

	\medspace
	
We note that, for $N\geq 3$, there holds $\frac{1}{4}\leq h_N^2(\lambda)\leq \left(\frac{N-2}{2} \right)^2$. Therefore, when $\lambda$ is in the range $[0,N-2),$ the term $(\frac{(N-2)^2}{4} - h^2_{N}(\lambda))$ is negative, but in this case we can combine the second term with third term  of R.H.S of \eqref{improved-poinc-lambda} and using the fact $\frac{1}{r^2} > \frac{1}{\sinh^2r},$ we  eventually obtain an  inequality which resembles (locally) Poincar\'e-Hardy inequality (See \cite[Remark~2.3]{BGGP}). This fact leads us a natural question whether we can improve \eqref{improved-poinc-lambda} for $\lambda \in [0,N-2).$ In Theorem~\ref{subspace_hardy} we prove this fact by choosing suitable subspaces of the Sobolev spaces.

\medskip 

As an immediate consequence of Theorem~\ref{subspace_hardy}, we establish an improved Hardy inequality on a certain subspace of the Sobolev spaces, namely $\mathcal{H}_0(\hn\setminus\{x_0\})$ (refer Section \ref{Section-2} for details) which is as follows: for all $N\geq 2$ there holds
\begin{align*}
	\int_{\hn}|\nabla_{\hn} u|^2\dv\geq\frac{N^2}{4}\int_{\hn}\frac{u^2}{r^2}\dv \text{ for all }u\in C_c^\infty(\mathcal{H}_0(\hn\setminus\{x_0\})).
\end{align*}
 Secondly, we prove a new kind of improved Caffarelli-Kohn-Nirenberg inequality in the hyperbolic space. 

\medspace

The paper is organized as follows: in Section \ref{Section-2} we introduce some of the notations and some 
geometric definitions related to the hyperbolic spaces and state our main Theorems and their consequences. Section \ref{section-3} is devoted to the proof of the results stated in Section \ref{Section-2}. Finally, in Section \ref{section-4} we prove the Caffarelli-Kohn-Nirenberg inequality on hyperbolic space involving radial derivatives.


\section{Preliminaries and Statement of main Results}\label{Section-2}

\subsection{Notations} It is well known that $N$-dimensional hyperbolic space $\hn$ admits polar coordinate decomposition structure. Namely we can write for some $x\in\hn$ as $x=(r, \theta_{1},\ldots, \theta_{N-1})\in(0,\infty)\times\mathbb{S}^{N-1}$, where $r$ denotes the geodesic distance between the point $x$ and a fixed pole $o$ in $\hn$ and $\mathbb{S}^{N-1}$ is the unit sphere in the $N$-dimensional euclidean space $\mathbb{R}^N$. 

Recall that the Riemannian Laplacian of a scalar function $u$ on $\hn$ is given by
\begin{align*}
	\Delta_{\hn} u (r, \theta_{1}, \ldots, \theta_{N-1})  &=
	\frac{1}{\sinh^2r} \frac{\partial}{\partial r} \left[ (\sinh r)^{N-1} \frac{\partial u}{\partial r}(r, \theta_{1},
	\ldots, \theta_{N-1}) \right]\\&
	+ \frac{1}{\sinh^2r} \Delta_{\mathbb{S}^{N-1}} u(r, \theta_{1}, \ldots, \theta_{N-1}),
\end{align*}
where $\Delta_{\mathbb{S}^{N-1}}$ is the Riemannian Laplacian on the unit sphere $\mathbb{S}^{N-1}$. In particular,
the radial contribution of the Riemannian Laplacian, namely the operator involving only on $r$, $\Delta_{r,\hn} u$ reads as
\begin{equation*}
	\Delta_{r,\hn} u = \frac{1}{(\sinh r)^{N-1}} \frac{\partial}{\partial r} \left[ (\sinh r)^{N-1} \frac{\partial u}{\partial r}
	 \right] =  u^{\prime \prime} + (N-1) (\coth r) u^{\prime},
\end{equation*}
where from now onwards {\it a prime \rm} will denote derivative w.r.t $r$. Also let us recall the Gradient in terms of the polar coordinate decomposition, given by
\begin{equation*}
	\nabla_{\hn}u(r, \theta_{1}, \ldots, \theta_{N-1})=\bigg(\frac{\partial u}{\partial r}(r, \theta_{1}, \ldots, \theta_{N-1}), \frac{1}{\sinh r}\nabla_{\mathbb{S}^{N-1}}u(r, \theta_{1}, \ldots, \theta_{N-1})\bigg),
\end{equation*}
where $\nabla_{\mathbb{S}^{N-1}}$ denotes the Gradient on the unit sphere $\mathbb{S}^{N-1}$. Again, the radial contribution of the Gradient, $\nabla_{r,\hn}u$ is defined as 
\begin{equation*}
	\nabla_{r,\hn}u=\bigg(\frac{\partial u}{\partial r}, 0 \bigg).
\end{equation*}

\medspace

\subsection{Spherical harmonics}\label{sph_hyp} Here we will recall some useful facts from spherical harmonics.
Let $u(x)=u(r,\sigma)\in C_{c}^\infty(\hn)$, $r\in({0},\infty)$ and $\sigma\in \mathbb{S}^{N-1}$, we can write
\begin{equation*}
	u(r,\sigma)=\sum_{n=0}^{\infty}a_n(r)P_n(\sigma)
\end{equation*}
in $L^2(\hn)$, where $\{ P_n \}$ is an orthonormal system of spherical harmonics and 
\begin{equation*}
	a_n(r)=\int_{\mathbb{S}^{N-1}}u(r,\sigma)P_n(\sigma) \ {\rm d}\sigma\,.
\end{equation*} 
A spherical harmonic $P_n$ of order $n$ is the restriction to $\mathbb{S}^{N-1}$ of a homogeneous harmonic polynomial of degree $n.$ Moreover, it
satisfies $$-\Delta_{\mathbb{S}^{N-1}}P_n(\sigma)=\lambda_n P_n(\sigma)$$
for all $n\in\mathbb{N}\cup\{0\}$, where $\lambda_n=(n^2+(N-2)n)$ are the eigenvalues of Laplace Beltrami operator $-\Delta_{\mathbb{S}^{N-1}}$ on $\mathbb{S}^{N-1}$ with corresponding eigenspace dimension $d_n$. We note that $\lambda_n\geq N-1$ for $n\geq 1$, $\lambda_0=0$, $d_0=1$, $d_1=N$ and for $n\geq 2$ 
\begin{equation*}
	d_n=\binom{N+n-1}{n}-\binom{N+n-3}{n-2}.
\end{equation*}

Using the above expression in terms of spherical harmonics we can write 
\begin{align*}
	|\nabla_{\hn}u|^2 =  \sum_{n = 0}^{\infty} {a_{n}^{\prime}}^2 P_{n}^2  + \frac{a_{n}^2}{\sinh^2r}  |\nabla_{\mathbb{S}^{N-1}}P_{n}|^2
\end{align*}
and
\begin{align*}
	(\Delta_{\hn} u)^2 & = \sum_{n = 0}^{\infty} \left( a_{n}^{\prime \prime} + (N-1) (\coth r) a_{n}^{\prime} \right)^2 P_{n}^2 + 
	\sum_{n = 0}^{\infty} \frac{a_{n}^2}{\sinh^4r} (\Delta_{\mathbb{S}^{N-1}} P_{n})^2 \\
	& + 2 \sum_{n = 0}^{\infty} \left( a_{n}^{\prime \prime} + (N-1) (\coth r) a_{n}^{\prime} \right) \frac{a_{n}}{\sinh^2r} (\Delta_{\mathbb{S}^{N-1}} P_{n}) P_{n}.\notag
\end{align*}

Along with this, the radial contribution of the operators will look as follows:
\begin{align*}
	|\nabla_{r,\hn}u|^2 =  \sum_{n = 0}^{\infty} {a_{n}^{\prime}}^2 P_{n}^2\quad	\text{ and } \quad(\Delta_{r,\hn} u)^2 & = \sum_{n = 0}^{\infty} \left( a_{n}^{\prime \prime} + (N-1) (\coth r) a_{n}^{\prime} \right)^2 P_{n}^2.
\end{align*}

\medspace

\subsection{Subspace}\label{sub} We have already introduced the spherical harmonics decomposition for a function $u\in C_c^\infty(\hn)$.  Now we will introduce some subspace of $C_c^\infty(\hn\setminus\{x_0\})$ and we will show that constants for most of the known inequalities can be larger than optimal one.  Here is the required subspace mention below. First consider $u\in C_c^\infty(\hn)$ and we can write $$u(x)=u(r,\sigma)=\sum_{n=0}^{\infty}a_n(r)P_n(\sigma)=\sum_{n=0}^{\infty}u_n(r,\sigma)$$
and now define for $j\in\mathbb{N}\cup \{0\}$ and $\Omega\subset\hn$,
\begin{align*}
	\mathcal{H}_j(\Omega)&=\biggl\{u\in C_c^\infty(\Omega) : u_0(r,\sigma)=u_1(r,\sigma)=\cdots=u_j(r,\sigma)=0\biggr\}\\
	&=\biggl\{u\in C_c^\infty(\Omega) : u(x)=u(r,\sigma)=\sum_{n=j+1}^{\infty}a_n(r)P_n(\sigma)\biggr\}.
\end{align*}
For notational clarity we want to mention that $\mathcal{H}_{-1}(\Omega)=C_c^\infty(\Omega).$

\medspace

Now we are in a situation to state our main results. 

\subsection{Improved Hardy-Rellich type inequality on subspace}\label{result_hyp} Our first aim of this paper is to show that the optimal constants of the Hardy type inequalities can indeed be improved on $\mathcal{H}_j(\hn\setminus\{x_0\}):$

\begin{theorem}\label{subspace_hardy}
	Let $N \geq 2$ and  $(r^{N-1}V,r^{N-1}W)$ be a Bessel pair on $(0,R)$ with $0 < R \leq \infty$ and $f$ be the corresponding positive solution. Then for all $u\in \mathcal{H}_j(B_R\setminus\{x_0\})$ there holds
	\begin{align*}
		&\int_{B_R}V(r)|\nabla_{\hn} u|^2\dv\geq \int_{B_R}W(r) \, u^2\dv+(j+1)(N+j-1)\int_{B_R}\frac{V(r) \, u^2}{\sinh^2 r}\dv\\& +\int_{B_R}V(r)(f(r))^2\bigg|\nabla_{r, \hn}\bigg(\frac{u}{f(r)}\bigg)\bigg|^2\dv-(N-1)\int_{B_R}V(r)\frac{f^\prime(r)}{f(r)}\bigg(\coth r - \frac{1}{r}\bigg) \, u^2 \, \dv.
	\end{align*}
\end{theorem}

\medspace

The above result can be compared with the version of \cite[Theorem 3.2]{FLLM} in hyperbolic space.  It is worth mentioning that if we replace the operator in terms of radial derivative, then our methods do not produce any improved version in the inequality and this suggests that the geometric derivative plays a crucial role here. This led us to the following theorem: 

\begin{theorem}\label{subspace_hardy_rem}
	Let $V(r)$ be a radial function on $(0,R)$. Then for all $u\in \mathcal{H}_j(B_R\setminus\{x_0\})$ there holds
	\begin{align*}
		&\int_{B_R}V(r)|\nabla_{\hn} u|^2\dv-\int_{B_R}V(r)\bigg(\frac{\partial u}{\partial r}\bigg)^2\dv= (j+1)(N+j-1)\int_{B_R}\frac{V(r)|u|^2}{\sinh^2 r}\dv.
	\end{align*}
\end{theorem}

\medspace

Using Theorem~\ref{subspace_hardy} we shall obtain the following corollary, one can compare the result below with \cite[Theorem 5.1]{EGR-21} and \cite[Theorem 2.1]{BGGP} :

\begin{corollary}\label{improved-hardy-radial}
	Let $N \geq 2.$ For all $0\leq  \lambda \leq  \lambda_{1}(\hn)$ and all $u\in \mathcal{H}_j(\hn\setminus\{x_0\})$ there holds
	\begin{align*}
		\int_{\hn}|\nabla_{\hn} u|^2\, \emph{d}v_{\hn} &\geq \lambda \int_{\hn}  |u|^2\, \emph{d}v_{\hn} +h_N^2(\lambda) \int_{\hn} \frac{|u|^2}{r^2} \ \emph{d}v_{\hn} \\
		& +  \left[\left(\frac{N}{2} \right)^2 - h_N^2(\lambda)+j(N+j)\right] \int_{\hn}  \frac{|u|^2}{\sinh^2 r} \ \emph{d}v_{\hn}\\&+  \gamma_{N}(\lambda)\, h_N(\lambda)\int_{\hn} \frac{r \coth r - 1}{r^2}\, |u|^2  \ \emph{d}v_{\hn} \notag \\
		&+\int_{\hn}(\Psi_{\lambda}(r))^2\bigg|\nabla_{\hn}\bigg(\frac{u}{\Psi_{\lambda}(r)}\bigg)\bigg|^2\dv, \notag
	\end{align*}
where $\gamma_{N}(\lambda):=\sqrt{(N-1)^2-4\lambda}$, $h_N(\lambda):=\frac{\gamma_{N}(\lambda)+1}{2}$ and $\Psi_{\lambda}(r) := r^{-\frac{N-2}{2}} \left(\frac{\sinh r}{r} \right)^{-\frac{N-1+\gamma_N(\lambda)}{2}}.$ 
\end{corollary}
\begin{remark}
	We can see that for all $\lambda\in[0,\lambda_{1}(\hn)]$, we have $$
	\frac{N^2}{4} - h_N^2(\lambda)+j(N+j)\geq 0, \quad \forall  \ j \geq 0. 
	$$
	Therefore the above inequality provides an improved constant than the Hardy constant  \cite[Theorem~2.1]{BGGP}. In particular for $u\in \mathcal{H}_0(\hn\setminus\{x_0\})$ and $\lambda=0$ there holds
	\begin{align*}
		\int_{\hn}|\nabla_{\hn} u|^2 \, \dv & \geq \frac{N^2}{4} \int_{\hn}\frac{|u|^2}{r^2} \, \dv.
	\end{align*}
\end{remark}

\medspace

Now substituting the  Bessel pair $\big(r^{N-1}r^{-\alpha},r^{N-1}\big(\frac{N-\alpha-2}{2}\big)^2r^{-\alpha-2}\big)$ in Theorem \ref{subspace_hardy} on the interval $(0,\infty)$ with positive solution $f(r)=r^{-\frac{N-\alpha-2}{2}}$ we obtain the following weighted Hardy inequality. 

\begin{corollary}\label{cor_wg_hardy}
	For all $u\in \mathcal{H}_j(B_R\setminus\{x_0\})$ there holds
	\begin{align*}
		\int_{\hn}\frac{|\nabla_{\hn} u|^2}{r^\alpha}\dv&\geq \bigg(\frac{N-\alpha-2}{2}\bigg)^2\int_{\hn}\frac{|u|^2}{r^{\alpha+2}}\dv\\&+(j+1)(N+j-1)\int_{\hn}\frac{|u|^2}{r^\alpha\sinh^2 r}\dv\\& +\int_{\hn}\frac{1}{r^{N-2}}\big|\nabla_{\hn}\big(r^{\frac{N-\alpha-2}{2}}u\big)\big|^2\dv\\&+\frac{(N-1)(N-\alpha-2)}{2}\int_{\hn}\frac{(r\coth r-1)}{r^{\alpha+2}}|u|^2\dv.
	\end{align*}
\end{corollary}

\medspace

\subsection{Caffarelli-Kohn-Nirenberg inequalities} This subsection is motivated by the recent progress in Caffarelli-Kohn-Nirenberg inequalities (shortly,  CKN inequalities),
see for instance (\cite{Caffa, CatWang, LamCazaFly, costa, dong, lamlu}).   The generalisation to  Riemannian manifolds, or more generally to metric space has been obtained by the authors in \cite{KO}, we also refer the works of  \cite{AKR, ST} which are closely related to CKN inequalities.  Very recently,  Nguyen \cite{VHNCKN} studied CKN inequalities in detail and proved several interesting results on manifolds with appropriate curvature assumptions. We obtain the following version of CKN inequality : 

\begin{theorem}\label{ckn}
	Let $N \geq 2,$  for all $u\in C_c^{\infty}(\hn\setminus \{x_0\})$, and $\alpha, \beta \in \mathbb{R},$ there holds
	\begin{align*}
		\bigg(\int_{\hn}\frac{|\nabla_{r,\hn}u|^2}{r^\alpha}\dv\bigg)&\bigg(\int_{\hn}r^{\alpha-2\beta+2}|u|^2\dv\bigg)\\&\geq\max\biggl\{\frac{(N-\beta)^2}{4},\frac{(N-2\alpha-\beta-4)^2}{4}\biggr\}\bigg(\int_{\hn}\frac{|u|^2}{r^{\beta}}\dv\bigg)^2.\nonumber
	\end{align*}
\end{theorem}

\medspace

To this end, we shall derive an improved version of Caffarelli-Kohn-Nirenberg inequality by providing an explicit remainder term.  The result can be stated as:
\begin{theorem}\label{ckn-rem}
	For all $u\in C_c^{\infty}(\hn\setminus \{x_0\})$, there holds
	\begin{align*}
		&\bigg(\int_{\hn}\frac{|\nabla_{r,\hn}u|^2}{r^\alpha}\dv\bigg)\bigg(\int_{\hn}r^{\alpha-2\beta+2}|u|^2\dv\bigg)-\frac{(N-\beta)^2}{4}\bigg(\int_{\hn}\frac{|u|^2}{r^{\beta}}\dv\bigg)^2\\&=\frac{\int_{\hn}\bigg|\bigg(\int_{\hn}\frac{|u|^2}{r^{2\beta-\alpha-2}}\dv\bigg)\frac{x}{r^{\frac{\alpha}{2}+1}}\frac{\partial u}{\partial r}-\bigg(\frac{N-\beta}{2}\int_{\hn}\frac{|u|^2}{r^{\beta}}\dv\bigg)\frac{x}{r^{\beta-\frac{\alpha}{2}}}u\bigg|^2\dv}{\int_{\hn}\frac{|u|^2}{r^{2\beta-\alpha-2}}\dv}.
	\end{align*}
\end{theorem}

\begin{remark}
The above inequalities should be compared with \cite[Theorem 1.3]{VHNCKN}, where under certain assumptions on $\alpha, \beta$ and dimension $N,$ a more general inequalities 
have been obtained but the above-stated inequalities hold without any restriction on dimensions and parameters $\alpha, \beta.$ Moreover a detailed discussion on the sharpness of the constant $\frac{(N-\beta)^2}{4}$ has been provided in \cite{VHNCKN}. 
\end{remark}

\subsection{Result holds for Riemannian symmetric manifold.} In this subsection, we briefly consider extensions of our results to 
 Riemannian symmetric model manifolds. One of the key methods exploited in this article is the \it spherical harmonic decomposition \rm 
 and in general, this is applicable to  \it, Riemannian models, \rm also. 

\medspace

An $N$-dimensional Riemannian model $\,(M,g,\psi)$ is an $N$-dimensional Riemannian manifold admitting a pole  $o\in M$ and whose metric $g$ in spherical coordinates around $o$ looks like
\begin{equation*}
	{\rm d}s^2 = {\rm d}r^2 + \psi^2(r) \, {\rm d}\omega^2,
\end{equation*}
where the coordinate $r,$ denotes the Riemannian distance from the fixed pole $o$ and ${\rm d}\omega^2$ represents the metric on the $N$-dimensional unit sphere $\mathbb{S}^{N-1}$ and the non-negative function $\psi$ satisfies
\begin{align*}
	\psi\in C^\infty([0,+\infty)), \,  \psi > 0 \text{ on } (0,+\infty),  \psi'(0) =1 \text{ and } \psi^{(2k)}(0)= 0 \text{ for all } k\geq 0\,.
\end{align*}
These conditions on $\psi$ confirm that the manifold is smooth and
the metric at the pole $o$ is given by the Euclidean metric. For example, all the assumptions above are satisfied 
by $\psi(r)=r$ and by $\psi(r)=\sinh(r)$: In the first case $M$ coincides with the Euclidean space $\mathbb{R}^N$, in the latter with the hyperbolic space $\hn$.

\medskip 

We emphasize that our arguments rely on the careful analysis of the polar coordinate representation of the Laplace-Beltrami and the Gradient operator on the hyperbolic space and exploiting the spectral analysis of $-\Delta_{\mathbb{S}^{N-1}}$. These fundamental operators on Riemannian models are given by: 
\begin{align*}
	\Delta_{g} \, = \,  \frac{\partial^2}{\partial \, r^2} \, + \, (N-1)\frac{\psi^\prime(r)}{\psi(r)}\, \frac{\partial}{\partial \, r} \, + \, \frac{1}{\psi^2(r)} \, \Delta_{\mathbb{S}^{N-1}}\, \text{ and }	\,\nabla_{g} \, = \,  \bigg(\frac{\partial}{\partial \, r}\,, \frac{1}{\psi(r)}\nabla_{\mathbb{S}^{N-1}}\bigg).
\end{align*}
Also the radial version of these operator is defined by 
\begin{align*}
	\Delta_{r,g} \, = \,  \frac{\partial^2}{\partial \, r^2} \, + \, (N-1)\frac{\psi^\prime(r)}{\psi(r)}\, \frac{\partial}{\partial \, r} \, \text{ and }	\,\nabla_{r,g} \, = \,  \bigg(\frac{\partial}{\partial \, r}\,,\,0\bigg).
\end{align*}

\medspace

It is worth noting that the method exploited in the Subsection~\ref{sph_hyp} can be extended to  Riemannian model $(M,g,\psi).$ Therefore, in the same way as described in the Subsection \ref{sub}, we can define the subspace $\mathcal{H}_j(B_R\setminus\{o\})$ of $C_c^\infty(B_R\setminus\{o\})$ for the model manifold $(M,g,\psi)$ with $B_R$ as an open ball of radius $R$ centred at pole $o$ of $M$. Hence, similar results like Theorem \ref{subspace_hardy} and \ref{subspace_hardy_rem} hold true in the model manifold involving the radial functions  $\psi(r)$ and $\psi^{\prime}(r)$. In particular,  the results describe as follows:

\begin{theorem}
	Let $(M,g,\psi)$ be a $N$-dimensional Riemannian model with $N\geq 2$ and $(r^{N-1}V,r^{N-1}W)$ be a Bessel pair on $(0,R)$ with $0 < R \leq \infty$ and $f$ be the corresponding positive solution. Then for all $u\in \mathcal{H}_j(B_R\setminus\{x_0\})$ there holds
	\begin{align*}
		&\int_{B_R}V(r)|\nabla_{g} u|^2\dvg\geq \int_{B_R}W(r) \, |u|^2\dvg+(j+1)(N+j-1)\int_{B_R}\frac{V(r) \, |u|^2}{\psi^2(r)}\dvg\\& +\int_{B_R}V(r)(f(r))^2\bigg|\nabla_{r, g}\bigg(\frac{u}{f(r)}\bigg)\bigg|^2\dvg-(N-1)\int_{B_R}V(r)\frac{f^\prime(r)}{f(r)}\bigg(\frac{\psi^\prime(r)}{\psi(r)} - \frac{1}{r}\bigg) \, |u|^2 \, \dvg.
	\end{align*}
\end{theorem}

\begin{theorem}
	Let $(M,g,\psi)$ be a $N$-dimensional Riemannian model with $N\geq 2$ and $V(r)$ be a radial function on $(0,R)$. Then for all $u\in \mathcal{H}_j(B_R\setminus\{x_0\})$ there holds
	\begin{align*}
		&\int_{B_R}V(r)|\nabla_{g} u|^2\dvg-\int_{B_R}V(r)\bigg(\frac{\partial u}{\partial r}\bigg)^2\dvg= (j+1)(N+j-1)\int_{B_R}\frac{V(r)|u|^2}{\psi^2(r)}\dvg.
	\end{align*}
\end{theorem}

\medspace

\begin{remark}

We shall not present the proof here as it can be easily done by choosing $\psi(r)$ instead of $\sinh(r)$ and exploiting \cite[Theorem 1.1]{FLLM} and following the method as in the case of hyperbolic space. The other derived remarks in the Subsection \ref{result_hyp} also hold in Riemannian models with obvious modification with appropriate assumptions on curvature.  For the sake of brevity, we shall skip the details. On the hand,  one of the main limitations of our method is it can not be generalised to general manifolds which do not admit symmetric structure. 
\end{remark}

\section{Proofs of Theorems \ref{subspace_hardy}, \ref{subspace_hardy_rem} and Corollary \ref{improved-hardy-radial}}\label{section-3}
{\bf Proof of Theorem \ref{subspace_hardy}.} 
	Start with $u\in \mathcal{H}_j(B_R\setminus\{x_0\})$. Then using spherical decomposition we have $$u(x)=u(r,\sigma)=\sum_{n=j+1}^{\infty}a_n(r)P_n(\sigma).$$ 
	Applying \cite[Theorem 3.2]{FLLM}, for each radial function $\{a_n\}$ we obtain
	\begin{align*}
		\int_{B_R}V(r)|\nabla_{\hn} a_n|^2\dv&= \int_{B_R}W(r)|a_n|^2\dv+\int_{B_R}V(r)(f(r))^2\bigg|\nabla_{\hn}\bigg(\frac{a_n}{f(r)}\bigg)\bigg|^2\dv\\&-(N-1)\int_{B_R}V(r)\frac{f^\prime}{f}\bigg(\coth r - \frac{1}{r}\bigg)|a_n|^2\dv.
	\end{align*}
	
	Using spherical harmonics and the fact that for  $n\geq j+1$ we have $\lambda_{n}\geq (j+1)(N+j-1),$ we can write 
	\begin{align*}
		&\int_{B_R}V(r)|\nabla_{\hn} u|^2\dv =\sum_{n =j+1}^{\infty}\bigg[\int_{0}^{R}V(r){a_n^\prime}^2(\sinh r)^{N-1}\dr+\lambda_{n}\int_{0}^{R}V(r)a_n^2(\sinh r)^{N-3}\dr\bigg]\\
		&= \sum_{n =j+1}^{\infty}\bigg[\int_{0}^{R}V(r){a_n^\prime}^2(\sinh r)^{N-1}\dr+ \lambda_n \int_{0}^{R}V(r)a_n^2(\sinh r)^{N-3}\dr\bigg]\\
		& = \sum_{n =j+1}^{\infty} \bigg[\int_{0}^{R} W(r) \, a_n^2(r) \, (\sinh r)^{N-1} \, {\rm d}r \, + \, \int_{0}^{R} V(r) \, (f(r))^2 \frac{\partial}{\partial r}\left( \frac{a_n}{f} \right)\, 
		(\sinh r)^{N-1} \, {\rm d}r  \bigg.  \\
		&  - (N-1) \int_{0}^{R} V(r) \, \frac{f^{\prime}}{f} \left( \coth r - \frac{1}{r} \right) \, a_n^2 \, (\sinh r)^{N-1} \, {\rm d}r 
		+ \lambda_n \int_{0}^{R}V(r)a_n^2(\sinh r)^{N-3}\dr \bigg] \\
		& \geq \int_{B_R}W(r)|u|^2\dv+(j+1)(N+j-1)\int_{B_R}\frac{V(r)|u|^2}{\sinh^2 r}\dv\\&+\int_{B_R}V(r)(f(r))^2\bigg|\nabla_{r, \hn}\bigg(\frac{u}{f(r)}\bigg)\bigg|^2\dv-(N-1)\int_{B_R}V(r)\frac{f^\prime}{f}\bigg(\coth r - \frac{1}{r}\bigg)|u|^2\dv.
	\end{align*} 

In the last line we have used the orthonormal properties of $\{P_n\}$ i.e. $\int_{\mathbb{S}^{N-1}}P_n(\sigma)P_m(\sigma)${\rm d}$\sigma=\delta_{nm}$ and writing back in terms of of $u$, we finish the prove. Precisely we have used the following identities:
\begin{align*}
\int_{B_R}W(r)|u|^2\dv=	\sum_{n =j+1}^{\infty} \int_{0}^{R} W(r) \, a_n^2(r) \, (\sinh r)^{N-1} \, {\rm d}r,
\end{align*}
\begin{align*}
	\int_{B_R}\frac{V(r)|u|^2}{\sinh^2 r}\dv=\sum_{n =j+1}^{\infty} \int_{0}^{R}V(r)a_n^2(\sinh r)^{N-3}\dr,
\end{align*}
\begin{align*}
	\int_{B_R}V(r)\frac{f^\prime}{f}\bigg(\coth r - \frac{1}{r}\bigg)|u|^2\dv=\sum_{n =j+1}^{\infty}\int_{0}^{R} V(r) \, \frac{f^{\prime}}{f} \left( \coth r - \frac{1}{r} \right) \, a_n^2 \, (\sinh r)^{N-1} \, {\rm d}r,
\end{align*}
and 
\begin{align*}
	\int_{B_R}V(r)(f(r))^2\bigg|\nabla_{r, \hn}\bigg(\frac{u}{f(r)}\bigg)\bigg|^2\dv=\sum_{n =j+1}^{\infty} \int_{0}^{R} V(r) \, (f(r))^2 \frac{\partial}{\partial r}\left( \frac{a_n}{f} \right)\, 
	(\sinh r)^{N-1} \, {\rm d}r.
\end{align*}
 \hfill $\Box$
 
 \medspace
 
 {\bf Proof of Theorem \ref{subspace_hardy_rem}.} The primary tool used in the proof is again the spherical decomposition technique. So begin with $u\in \mathcal{H}_j(B_R\setminus\{x_0\})$ and we can write $$u(x)=u(r,\sigma)=\sum_{n=j+1}^{\infty}a_n(r)P_n(\sigma).$$ 
 Now we will break every term of the required identity individually in terms of spherical harmonics and after simple modification, we have the following identities:
 \begin{align*}
 	\int_{B_R}V(r)|\nabla_{\hn} u|^2\dv=\sum_{n =j+1}^{\infty}\bigg[\int_{0}^{R}V(r){a_{n}^{\prime}}^2(\sinh r)^{N-1}\dr+\lambda_{n}\int_{0}^{R}V(r)a_n^2(\sinh r)^{N-3}\dr\bigg],
 \end{align*}
\begin{align*}
	\int_{B_R}V(r)\bigg(\frac{\partial u}{\partial r}\bigg)^2\dv=\sum_{n =j+1}^{\infty}\int_{0}^{R}V(r){a_{n}^{\prime}}^2(\sinh r)^{N-1}\dr,
\end{align*}
and
\begin{align*}
	\int_{B_R}\frac{V(r)|u|^2}{\sinh^2 r}\dv=\sum_{n =j+1}^{\infty} \int_{0}^{R}V(r)a_n^2(\sinh r)^{N-3}\dr.
\end{align*}

\medspace
 
Therefore using these above identities on the underlying subspace we shall obtain
	\begin{align*}
		&\int_{B_R}V(r)|\nabla_{\hn} u|^2\dv-\int_{B_R}V(r)\bigg(\frac{\partial u}{\partial r}\bigg)^2\dv=\sum_{n =j+1}^{\infty}\lambda_{n}\int_{0}^{R}V(r)a_n^2(\sinh r)^{N-3}\dr\\&
		= (j+1)(N+j-1)\int_{B_R}\frac{V(r)|u|^2}{\sinh^2 r}\dv+\sum_{n =j+1}^{\infty}(\lambda_{n}-\lambda_{j+1})\int_{0}^{R}V(r)a_n^2(\sinh r)^{N-3}\dr\\&= (j+1)(N+j-1)\int_{B_R}\frac{V(r)|u|^2}{\sinh^2 r}\dv+\sum_{n =j+1}^{\infty}(n-j-1)(N+n+j-1)\int_{B_R}\frac{V(r)|u|^2}{\sinh^2 r}\dv.
	\end{align*}
	In the end substituting the value of $\lambda_{n}=n(N+n-2)$ we conclude the proof. \hfill $\Box$
	
\medspace
	
{\bf Proof of Corollary \ref{improved-hardy-radial}.}	The proof follows by substituting particular choice of family of Bessel pairs in Theorem \ref{subspace_hardy}. We recall a new kind 
of Bessel pair which was first introduced in \cite{EGR-21}, i.e.,  Bessel pairs $(r^{N-1},r^{N-1}W_{\lambda})$ with $0\leq  \lambda \leq  \lambda_{1}(\hn)$ and
\begin{align}\label{potential}
W_{\lambda}(r)&:=\lambda+h_N^2(\lambda) \frac{1}{r^2}+ \left(\left(\frac{N-2}{2} \right)^2-h_N^2(\lambda)\right)   \frac{1}{\sinh^2 r}\notag \\&+ \left( \frac{\gamma_{N}(\lambda)\,h_N(\lambda)}{r}+(N-1) \frac{\Psi_{\lambda}'(r)}{\Psi_{\lambda}(r)}\right) \bigg(\coth r - \frac{1}{r}\bigg) \qquad (r>0)\,,
\end{align}
where $ \gamma_{N}(\lambda)$ and $h_N(\lambda)$ are as defined in the Corollary of Theorem \ref{improved-hardy-radial} and
$$
\Psi_{\lambda}(r) := r^{-\frac{N-2}{2}} \left(\frac{\sinh r}{r} \right)^{-\frac{N-1+\gamma_N(\lambda)}{2}}\qquad (r>0)\,.
$$
One can see after a straightforward computation 
\begin{align*}
\Psi''_{\lambda}(r) &= \Psi_{\lambda}(r)\Big[ \frac{(1-N-\gamma_N(\lambda))^2}{4}+\frac{\gamma_N^2(\lambda)-1}{r^2}\\
&-\frac{(1-N-\gamma_N(\lambda))(1+N+\gamma_N(\lambda))}{4 \sinh^2 r} + \frac{(1-N-\gamma_N(\lambda)) h_N(\lambda) \coth r}{r} \Big]
\end{align*}
and recalling the definition of $\gamma_N(\lambda)$, it follows that $\Psi_{\lambda}(r)$ satisfies
$$(r^{N-1} \Psi_{\lambda}'(r))'+r^{N-1}W_{\lambda}(r) \Psi_{\lambda}(r)=0 \quad \text{for } r>0\,,$$
 namely $(r^{N-1},r^{N-1}W_{\lambda})$ is a Bessel pair with positive solution $\Psi_{\lambda}(r)$.  \hfill $\Box$


\section{Proof of Theorem~\ref{ckn} and Theorem~\ref{ckn-rem}}\label{section-4}

This section is devoted to the proof of Theorem~\ref{ckn}. The proof rests on the geometry of the hyperbolic space. 

\medspace 	

{\bf Proof of Theorem \ref{ckn}.} It is easy to see that  
	\begin{equation*}
		\frac{\partial u}{\partial r} = \frac{x}{r}\cdot \nabla_{\hn} u \, \, \text{ with }\, \, r=\varrho(x,x_0).
	\end{equation*}
	Hence we can write 
	\begin{equation*}
		|\nabla_{r, \hn}u|=\frac{|x\cdot\nabla_{\hn}u|}{r}.
	\end{equation*}
	Now start with $\mathcal{A}:=\int_{\hn}\big|r^{-\frac{\alpha}{2}-1}x\cdot\nabla_{\hn}u+tr^{-\beta+\frac{\alpha}{2}+1}u+sr^{-\frac{\alpha}{2}-1}u\big|^2\dv$ non-negative number for all real $s,t$ (to be chosen later). Here expanding $\mathcal{A}$, we deduce
	\begin{align*}
		\mathcal{A}&=\int_{\hn}\bigg[\frac{|\nabla_{r,\hn}u|^2}{r^\alpha}+t^2r^{\alpha-2\beta+2}|u|^2+s^2r^{-\alpha-2}|u|^2\bigg]\dv\\&+2t\int_{\hn}\frac{u(x\cdot\nabla_{\hn}u)}{r^{\beta}}\dv+2s\int_{\hn}\frac{u(x\cdot\nabla_{\hn}u)}{r^{\alpha+2}}\dv+2st\int_{\hn}\frac{|u|^2}{r^{\beta}}\dv.
	\end{align*}
	Now using divergence theorem we have for any real number $m$
	\begin{align*}
		\int_{\hn}\frac{u(x\cdot\nabla_{\hn}u)}{r^{m}}\dv=-\frac{(N-m)}{2}\int_{\hn}\frac{|u|^2}{r^{m}}\dv.
	\end{align*}
	Exploiting these into above and further computing  we deduce
	\begin{align*}
		\mathcal{A}&=\int_{\hn}\frac{|\nabla_{r,\hn}u|^2}{r^\alpha}\dv+t^2\int_{\hn}r^{\alpha-2\beta+2}|u|^2\dv+t(2s-N+\beta)\int_{\hn}\frac{|u|^2}{r^{\beta}}\dv\\&+s(s-N+\alpha+2)\int_{\hn}\frac{|u|^2}{r^{\alpha+2}}\dv.
	\end{align*}
	For our result we want last term to be zero and so two cases arise.\\
	In the case $s=0$ we deduce 
	\begin{align*}
		\mathcal{A}&=t^2\int_{\hn}r^{\alpha-2\beta+2}|u|^2\dv-t(N-\beta)\int_{\hn}\frac{|u|^2}{r^{\beta}}\dv+\int_{\hn}\frac{|\nabla_{r,\hn}u|^2}{r^\alpha}\dv.
	\end{align*}
	The above quadratic equation in terms of $t$ is always nonnegative for every real number $t$ hence it's discriminant will be non-positive and we deduce
	\begin{align*}
		\bigg(\int_{\hn}\frac{|\nabla_{r,\hn}u|^2}{r^\alpha}\dv\bigg)\bigg(\int_{\hn}r^{\alpha-2\beta+2}|u|^2\dv\bigg)\geq\frac{(N-\beta)^2}{4}\bigg(\int_{\hn}\frac{|u|^2}{r^{\beta}}\dv\bigg)^2.
	\end{align*}
	For the second case $s=N-\alpha-2$ we deduce 
	\begin{align*}
		\mathcal{A}&=t^2\int_{\hn}r^{\alpha-2\beta+2}|u|^2\dv+t(N-2\alpha-\beta-4)\int_{\hn}\frac{|u|^2}{r^{\beta}}\dv+\int_{\hn}\frac{|\nabla_{r,\hn}u|^2}{r^\alpha}\dv.
	\end{align*}
	The above quadratic equation in terms of $t$ is always nonnegative for every real number $t$ hence it's discriminant will be non-positive and we deduce
	\begin{align*}
		\bigg(\int_{\hn}\frac{|\nabla_{r,\hn}u|^2}{r^\alpha}\dv\bigg)\bigg(\int_{\hn}r^{\alpha-2\beta+2}|u|^2\dv\bigg)\geq\frac{(N-2\alpha-\beta-4)^2}{4}\bigg(\int_{\hn}\frac{|u|^2}{r^{\beta}}\dv\bigg)^2.
	\end{align*}
	Finally combining these two cases we got the result.   \hfill $\Box$
	
	\medspace
	
	{\bf Proof of Theorem \ref{ckn-rem}.} We will start with for every $t\in\mathbb{R}$, there holds
	\begin{align*}
		&\int_{\hn}\bigg|\frac{x}{r^{\frac{\alpha}{2}+1}}\frac{\partial u}{\partial r}+t\frac{x}{r^{\beta-\frac{\alpha}{2}}}u\bigg|^2\dv\\&=\int_{\hn}\frac{1}{r^{\alpha}}\bigg(\frac{\partial u}{\partial r}\bigg)^2\dv+t^2\int_{\hn}\frac{|u|^2}{r^{2\beta-\alpha-2}}\dv+2t\int_{\hn}\frac{u}{r^{\beta-1}}\frac{\partial u}{\partial r}\dv.
	\end{align*}

Now using divergence theorem we have 
\begin{align*}
	\int_{\hn}\frac{u}{r^{\beta-1}}\frac{\partial u}{\partial r}\dv=\int_{\hn}\frac{u(x\cdot\nabla_{\hn}u)}{r^{\beta}}\dv=-\frac{(N-\beta)}{2}\int_{\hn}\frac{|u|^2}{r^{\beta}}\dv.
\end{align*}

By choosing 
$t:=\frac{(N-\beta)}{2}\frac{\int_{\hn}\frac{|u|^2}{r^{\beta}}\dv}{\int_{\hn}\frac{|u|^2}{r^{2\beta-\alpha-2}}\dv},
$
 we deduce
\begin{align*}
&\int_{\hn}\frac{|\nabla_{r,\hn}u|^2}{r^\alpha}\dv+	\frac{(N-\beta)^2}{4}\frac{(\int_{\hn}\frac{|u|^2}{r^{\beta}}\dv)^2}{(\int_{\hn}\frac{|u|^2}{r^{2\beta-\alpha-2}}\dv)^2}\bigg(\int_{\hn}\frac{|u|^2}{r^{2\beta-\alpha-2}}\dv\bigg)\\&-\frac{(N-\beta)^2}{2}\frac{(\int_{\hn}\frac{|u|^2}{r^{\beta}}\dv)^2}{\int_{\hn}\frac{|u|^2}{r^{2\beta-\alpha-2}}\dv}=\int_{\hn}\bigg|\frac{x}{r^{\frac{\alpha}{2}+1}}\frac{\partial u}{\partial r}+t\frac{x}{r^{\beta-\frac{\alpha}{2}}}u\bigg|^2\dv,
\end{align*}
which in turn implies,
\begin{align*}
& \bigg(\int_{\hn}\frac{|\nabla_{r,\hn}u|^2}{r^\alpha}\dv\bigg)\bigg(\int_{\hn}\frac{|u|^2}{r^{2\beta-\alpha-2}}\dv\bigg)-\frac{(N-\beta)^2}{4}\bigg(\int_{\hn}\frac{|u|^2}{r^{\beta}}\dv\bigg)^2\\& = \bigg(\int_{\hn}\frac{|u|^2}{r^{2\beta-\alpha-2}}\dv\bigg) \int_{\hn}\bigg|\frac{x}{r^{\frac{\alpha}{2}+1}}\frac{\partial u}{\partial r}+\frac{(N-\beta)}{2}\frac{\int_{\hn}\frac{|u|^2}{r^{\beta}}\dv}{\int_{\hn}\frac{|u|^2}{r^{2\beta-\alpha-2}}\dv}\frac{x}{r^{\beta-\frac{\alpha}{2}}}u\bigg|^2\dv,
\end{align*} 

and after rearranging terms in the last line we obtain the desired result. \hfill $\Box$
	
\medspace
   \par\bigskip\noindent
\textbf{Acknowledgments.}
D.~Ganguly is partially supported by the INSPIRE faculty fellowship (IFA17-MA98). P.~Roychowdhury is supported by the SRF fellowship (IISER-P/Proj./Offer/061/2022) under the funded project (Code: 31821586).

\end{document}